\documentclass[reqno]{amsart}
\usepackage{amsmath}

\def\eqn#1{(\ref{eq:#1})}

\newcommand{\qBin}[3]{\genfrac{[}{]}{0pt}{0}{#1}{#2}_{#3}}

\newcommand{\qBinpwr}[4]{\genfrac{[}{]}{0pt}{0}{#1}{#2}_{#3}^{#4}}

\theoremstyle{plain}
\newtheorem{thm}{Theorem}
\theoremstyle{remark}
\newtheorem*{rem}{Remark}

\numberwithin{equation}{section}
\allowdisplaybreaks[2]

\begin{document}

\title[The WP -- Bailey tree and its implications]
      {The WP -- Bailey tree and its implications}
\author[George~E.~Andrews]{George Andrews}
\address{Department of Mathematics,
         The Pennsylvania State University,
         University Park, PA~16802}
\email{andrews@math.psu.edu}
\author[Alexander~Berkovich]{Alexander Berkovich}
\address{Department of Mathematics, 
         University of Florida,
         Gainesville, FL~32611}
\email{alexb@math.ufl.edu}


\begin{abstract}
Our object is a thorough analysis of the WP-Bailey tree, a recent extension 
of classical Bailey chains. We begin by observing how the WP--Bailey tree 
naturally entails a finite number of classical $q$-hypergeometric 
transformation formulas. We then show how to move beyond this closed set of results 
and in the process we explicate heretofore mysterious identities of D.M.~Bressoud. 
Next, we use WP--Bailey pairs to provide a new proof of recent formula of A.N.~Kirillov. 
Finally, we discuss the relation between our approach and that of W.H.~Burge.
\end{abstract}

\maketitle

\section{Introduction}
\label{sec:1}

A classical Bailey pair \cite{An1} is a sequence $(\alpha_n,\beta_n)$ of pairs of 
rational functions of several complex variables subject to the identity
\begin{equation}
\beta_n=\sum_{j=0}^n\frac{\alpha_j}{(q;q)_{n-j}(aq;q)_{n+j}},
\label{eq:1.1}
\end{equation}
where
\begin{align}
(a;q)_n\equiv (a)_n=
\begin{cases} (1-a)(1-aq)\cdots(1-aq^{n-1}), & \mbox{if } n>0, \\
               1, & \mbox{if } n=0. \end{cases}
\label{eq:1.2}
\end{align}

\noindent
From such Bailey pairs $(\alpha_n,\beta_n)$ it is possible to construct new 
Bailey pairs by the discovery \cite{An1} that $(\alpha'_n,\beta'_n)$ are also 
Bailey pairs where
\begin{equation}
\alpha^\prime_n=\frac{(\rho_1,\rho_2)_n}{(\frac{aq}{\rho_1},\frac{aq}{\rho_2})_n}
(\frac{aq}{\rho_1\rho_2})^n\alpha_n
\label{eq:1.3}
\end{equation}
and
\begin{equation}
\beta^\prime_n=\frac{1}{(\frac{aq}{\rho_1},\frac{aq}{\rho_2})_n}
\sum_{j=0}^n(\rho_1,\rho_2)_j 
\frac{\left(\frac{aq}{\rho_1\rho_2}\right)_{n-j}}{(q)_{n-j}}
\left(\frac{aq}{\rho_1\rho_2}\right)^j\beta_j
\label{eq:1.4}
\end{equation}

\noindent
with
\begin{equation}
(a_1,a_2,\ldots,a_s;q)_j\equiv(a_1,a_2,\ldots,a_s)_j=(a_1)_j(a_2)_j\cdots(a_s)_j.
\label{eq:1.5}
\end{equation}

\noindent
This construct has wide applications in number theory, analysis, physics 
and has been the topic of two recent survey articles \cite{An2} and \cite{War1}.
The final sections of the former article were devoted to the further development 
of ideas nascent in the work of Bailey (\cite{Ba1}, \S 9). Namely, following 
Bressoud \cite{Br}, one may add a further parameter, say $k$, to the definitions of
Bailey pair to form what was termed WP--Bailey pair in \cite{An2}:
\begin{equation}
\beta_n(a,k)=\sum_{j=0}^n\frac{(\frac{k}{a})_{n-j}}{(q)_{n-j}}
\frac{(k)_{n+j}}{(aq)_{n+j}} \alpha_j(a,k).
\label{eq:1.6}
\end{equation}

\noindent
It was shown in \cite{An2}, that in this case, there are two distinct ways of 
constructing new WP-Bailey pairs. Both $(\alpha^\prime_n,\beta^\prime_n)$ and
$(\tilde\alpha_n,\tilde\beta_n)$ form a WP pair where
\begin{equation}
\alpha^\prime_n(a,k)=\frac{(\rho_1,\rho_2)_n}{(\frac{aq}{\rho_1},\frac{aq}{\rho_2})_n}
(\frac{k}{c})^n\alpha_n(a,c),
\label{eq:1.7}
\end{equation}

\begin{equation}
\beta^\prime_n(a,k)=\frac{(\frac{k\rho_1}{a},\frac{k\rho_2}{a})_n}
{(\frac{aq}{\rho_1},\frac{aq}{\rho_2})_n}\sum_{j=0}^n
\frac{(\rho_1,\rho_2)_j}{(\frac{k\rho_1}{a},\frac{k\rho_2}{a})_j}
\frac{1-cq^{2j}}{1-c} \frac{(\frac{k}{c})_{n-j}}{(q)_{n-j}}
\frac{(k)_{n+j}}{(qc)_{n+j}} (\frac{k}{c})^j \beta_j(a,c),
\label{eq:1.8}
\end{equation}
with 
\begin{equation}
c=\frac{k\rho_1\rho_2}{aq} 
\label{eq:1.9}
\end{equation}
while
\begin{equation}
\tilde\alpha_n(a,k)=\frac{(\frac{qa^2}{k})_{2n}}{(k)_{2n}}
\left(\frac{k^2}{qa^2}\right)^n \alpha_n\left(a,\frac{qa^2}{k}\right),
\label{eq:1.10}
\end{equation}

\begin{equation}
\tilde\beta_n(a,k)=\sum_{j=0}^n \frac{(\frac{k^2}{qa^2})_{n-j}}{(q)_{n-j}}
\left(\frac{k^2}{qa^2}\right)^j \beta_j\left(a,\frac{qa^2}{k}\right).
\label{eq:1.11}
\end{equation}

\begin{rem}
It is important to observe that with a double application of the construct \eqn{1.10} 
and \eqn{1.11} we return to the original WP--Bailey pair. Also, we note that 
\eqn{1.6}--\eqn{1.8} with $k=0$ give back \eqn{1.1}--\eqn{1.4}.
\end{rem}

In Section 2, we recount the essential results from \cite{An2}. The main point 
here is that while \cite{An2} brought a collection of $q$-hypergeometric identities 
under the classical Bailey umbrella, nonetheless this is only the beginning.

In Section 3, we undertake a careful study of how we may utilize WP pairs to generate 
further identities of the classical form wherein no multiple series appear. The surprise 
in this section is that the natural iterative process causes the sequence of results to
repeat. We conclude this section by presenting two new identities of the classical form.

Section 4 recaps the work of Bressoud in \cite{Br}, which contains  
surprising polynomial identities including particularly attractive polynomial 
refinements of the Rogers--Ramanujan identities. In Section 4, we identify Bressoud's 
Bailey Lemma with the limiting case of our construct \eqn{1.10}, \eqn{1.11} and comment 
on the three WP-pairs introduced in \cite{Br}.
Next, we explore the $q$-hypergeometric consequences of Bressoud's WP-pairs and
derive four more new $q$-hypergeometric transformations of the classical form.

In Section 5, we present a number of doubly bounded polynomial identities of the 
Rogers--Ramanujan type, making contacts with the work of Bressoud \cite{Br} and 
Warnaar \cite{War2}. 

Section 6 is devoted to an identity of Kirillov. We note that Kirillov mentioned 
this identity at a Special Functions 2000 Conference in Tempe. It was in an effort to 
understand Kirillov's discovery that the WP-Bailey tree was invented.

In Section 7, we briefly review the work of Burge and then, identify
certain WP-pairs as Burge pairs. This enables us to derive a number of identities
for $q$-multiple series. For background and recent work on multiple series Rogers--Ramanujan 
type identities reader may consult \cite{An1}, \cite{An3}, \cite{BM1}--\cite{BP2},  
\cite{Br2}, \cite{F}, \cite{War2}.

In Section 8, we modify the construct \eqn{1.10}, \eqn{1.11} in such a way 
that it can be applied to general Burge pairs.

Finally, Section 9 contains our concluding remarks.

Let us now recall some standard $q$-hypergeometric definitions and notations.\\
The generalized basic hypergeometric function is denoted by

\begin{equation}
_{r+1}\phi_r \left(\begin{array}{l} a_1,a_2,\ldots,a_{r+1};q,z \\
b_1,b_2,\ldots,b_r \end{array}\right)=
\sum_{j=0}^\infty \frac{(a_1,a_2,\ldots,a_{r+1})_j}{(q,b_1,b_2,\ldots,b_r)_j}z^j.
\label{eq:1.12}
\end{equation}

\noindent
We shall call a basic hypergeometric function well-poised if the parameters satisfy 
the relations

\begin{equation}
qa_1=b_1a_2=b_2a_3=\cdots=b_ra_{r+1},
\label{eq:1.13}
\end{equation}
and very-well-poised if, in addition,

\begin{equation}
a_2=q\sqrt a_1, a_3=-q\sqrt a.
\label{eq:1.14}
\end{equation}
A nearly-poised series of the first kind is one which satisfies 

\begin{equation}
qa_1\neq b_1a_2=b_2a_3=\cdots=b_ra_{r+1},
\label{eq:1.15}
\end{equation}
whereas a nearly-poised series of the second kind satisfies

\begin{equation}
qa_1=b_1a_2=b_2a_3=\cdots=b_{r-1}a_r\neq b_ra_{r+1}.
\label{eq:1.16}
\end{equation}

\noindent
In order to simplify some of the formulas involving very-well-poised 
$_{r+1}\phi_r$-series we shall frequently use the compact notation
\begin{equation}
_{r+1}W_r (a_1;a_4,a_5,\ldots,a_{r+1};q,z)=_{r+1}\phi_r
\left(\begin{array}{l} a_1,q\sqrt a_1,-q\sqrt a_1,a_4,a_5,\ldots,a_{r+1};q,z \\
\sqrt a_1,-\sqrt a_1,\frac{qa_1}{a_4},\frac{qa_1}{a_5},\cdots,\frac{qa_1}{a_{r+1}}
\end{array}\right).
\label{eq:1.17}
\end{equation}
Next, we introduce bibasic hypergeometric series

$$
\Phi\left[\begin{array}{l} a_1,a_2,\ldots,a_r:c_1,c_2,\ldots,c_s;q_1,q_2;z \\
b_1,b_2,\ldots,b_{r-1}:d_1,d_2,\ldots,d_s \end{array}\right]= 
$$
\begin{equation}
\sum_{j=0}^\infty 
\frac{(a_1,a_2,\ldots,a_r;q_1)_j}{(q_1,b_1,b_2,\ldots,b_{r-1};q_1)_j}
\frac{(c_1,c_2,\ldots,c_s;q_2)_j}{(d_1,d_2,\ldots,d_s;q_2)_j}z^j.
\label{eq:1.18}
\end{equation}

If $_{r+1}\phi_r$ series terminates, then by reversing the order of summation 
one can show that

$$
_{r+1}\phi_r\left(\begin{array}{l} a_1,\ldots,a_r,q^{-n};q,z \\
b_1,\ldots,b_r \end{array}\right)=\frac{(a_1,\ldots,a_r)_n}{(b_1,\ldots,b_r)_n}
(-1)^n q^{\frac{n(1-n)}{2}}
$$
\begin{equation}
\times_{r+1}\phi_r\left(\begin{array}{l} \frac{q^{1-n}}{b_1},\cdots,\frac{q^{1-n}}{b_r},
q^{-n};q,\frac{b_1,\ldots,b_r}{a_1,\ldots,a_r}\frac{q^{n+1}}{z} \\
\frac{q^{1-n}}{a_1},\ldots,\frac{q^{1-n}}{a_r} \end{array}\right)
\label{eq:1.19}
\end{equation}
with $n\in\Bbb{Z}_{\geq 0}$.
Using \eqn{1.19} it is easy to see that a terminating nearly-poised series of 
the second kind can be expressed as a multiple of a nearly-poised series of 
the first kind.

Finally, we remark that the acronym ``WP" stands for ``well poised". Indeed, let us
rewrite \eqn{1.16} as

\begin{equation}
\beta_n(a,k)=\frac{(k,\frac{k}{a})_n}{(q,aq)_n}\sum_{j=0}^n
\frac{(q^{-n},kq^n)_j}{(aq^{1+n},\frac{a}{k}q^{1-n})_j}
\left(\frac{aq}{k}\right)^j\alpha_j(a,k).
\label{eq:1.20}
\end{equation}
Clearly, in the above
\begin{equation}
(q^{-n})(aq^{1+n})= (kq^n)\left(\frac{a}{k}q^{1-n}\right)=aq.
\label{eq:1.21}
\end{equation}

\medskip

\section{The WP--Bailey tree}
\label{sec:2}

In the introduction, equation \eqn{1.6} provides a fundamental definition of a
WP--Bailey pair.
In \cite{An2}, it is shown (using only mathematical induction) that $(\alpha_n,\delta_{n,0})$ 
form a WP-Bailey pair, where

\begin{align}
\delta_{n,0}=
\begin{cases} 1, & \mbox{if } n=0, \\
              0, & \mbox{otherwise }, \end{cases}
\label{eq:2.1}
\end{align}
and
\begin{equation}
\alpha_n(a,k)=\frac{(a)_n(1-aq^{2n})}{(q)_n(1-a)}\frac{\left(\frac{a}{k}\right)_n}{(kq)_n}
\left(\frac{k}{a}\right)^n.
\label{eq:2.2}
\end{equation}
Applying \eqn{1.7} and \eqn{1.8} to produce a new WP-pair $(\alpha'_n,\beta'_n)$, we find

\begin{equation}
\alpha^\prime_n=\frac{(a,q\sqrt a,-q\sqrt a,\rho_1,\rho_2,\frac{a}{c})_n} 
{(q,\sqrt a,-\sqrt a,\frac{aq}{\rho_1},\frac{aq}{\rho_2},qc)_n}
\left(\frac{k}{a}\right)^n,
\label{eq:2.3}
\end{equation}
\begin{equation}
\beta^\prime_n=\frac{\left(\frac{k\rho_1}{a},\frac{k\rho_2}{a},k,\frac{k}{c}\right)_n}
{\left(\frac{aq}{\rho_1},\frac{aq}{\rho_2},q,qc\right)_n}
\label{eq:2.4}
\end{equation}
with $c$ as in \eqn{1.9}.
\begin{rem}
The WP--Bailey pair $(\alpha'_n,\beta'_n)$ was first discovered in the work of Singh \cite{Si}.
\end{rem}
If we substitute this latter pair into \eqn{1.6}, the resulting identity \cite{An2} is a 
well-known $q$-analog of Dougall's summation formula due to Jackson:
\begin{equation}
_8W_7(a;\rho_1,\rho_2,\frac{a}{c},kq^n,q^{-n};q,q)=
\frac{\left(aq,\frac{k}{c},\frac{k\rho_1}{a},\frac{k\rho_2}{a}\right)_n}
     {\left(cq,\frac{k}{a},\frac{aq}{\rho_1},\frac{aq}{\rho_2}\right)_n}.
\label{eq:2.5}
\end{equation}
Continuing the tree, we apply \eqn{1.7} and \eqn{1.8} to the WP-pairs 
\eqn{2.3} and \eqn{2.4} to obtain
\begin{equation}
\alpha^{\prime\prime}_n(a,k)=
\frac{(a,q\sqrt a,-q\sqrt a,\sigma_1,\sigma_2,\rho_1,\rho_2,\frac{ak}{c\tilde c})_n} 
{(q,\sqrt a,-\sqrt a,\frac{aq}{\sigma_1},\frac{aq}{\sigma_2},
\frac{aq}{\rho_1},\frac{aq}{\rho_2},\frac{qc\tilde c}{k})_n}
\left(\frac{k}{a}\right)^n,
\label{eq:2.6}
\end{equation}
$$
\beta^{\prime\prime}_n(a,k)=\frac{\left(\frac{k\sigma_1}{a},\frac{k\sigma_2}{a}\right)_n}
{\left(\frac{aq}{\sigma_1},\frac{aq}{\sigma_2}\right)_n}
\sum_{j=0}^n\frac{\left(\frac{k}{\tilde c}\right)_{n-j}}{(q)_{n-j}}
\frac{(k)_{n-j}}{(\tilde cq)_{n-j}}\frac{1-\tilde cq^{2j}}{1-\tilde c}
$$
\begin{equation}
\times\frac{\left(\sigma_1,\sigma_2,\frac{\tilde c\rho_1}{a},\frac{\tilde c\rho_2}{a},
\tilde c,\frac{k}{c}\right)_j}{\left(
\frac{k\sigma_1}{a},\frac{k\sigma_2}{a},\frac{aq}{\rho_1},\frac{aq}{\rho_2},
\frac{qc\tilde c}{k},q\right)_j}\left(\frac{k}{\tilde c}\right)^j,
\label{eq:2.7}
\end{equation}
where $\sigma_1,\sigma_2$ are new free parameters, $c$ is defined in \eqn{1.9}, and
\begin{equation}
\tilde c=\frac{k\sigma_1\sigma_2}{aq}.
\label{eq:2.8}
\end{equation}
We now substitute the above WP-pair into \eqn{1.6} to derive \cite{An2}:
$$
_{10}W_9(a;\rho_1,\rho_2,\frac{ak}{c\tilde c},\sigma_1,\sigma_2,kq^n,q^{-n};q,q)=
\frac{\left(aq,\frac{k}{\tilde c},\frac{k\sigma_1}{a},\frac{k\sigma_2}{a}\right)_n}
     {\left(\tilde cq,\frac{k}{a},\frac{aq}{\sigma_1},\frac{aq}{\sigma_2}\right)_n}
$$
\begin{equation}
\times_{10}W_9\left(\tilde c;\frac{\tilde c\rho_1}{a},\frac{\tilde c\rho_2}{a},\frac{k}{c},
\sigma_1,\sigma_2,kq^n,q^{-n};q,q\right).
\label{eq:2.9}
\end{equation}
It is easy to recognize formula \eqn{2.9} as Bailey's imposing 
$_{10}W_9\rightarrow\mbox{}_{10}W_9$ transformation. This result of Bailey 
has turned out to be one of the most fecund in the entire subject. 
Indeed, Gasper and Rahman \cite{GR} devote nearly five pages to its implications.

Finally, in our recap of the results sketched in \cite{An2}, we note that if we apply 
\eqn{1.10} and \eqn{1.11} using the WP-pair \eqn{2.3}, \eqn{2.4}, we obtain
\begin{equation}
\tilde\alpha_n(a,k)=\frac{(a)_n}{(q)_n}\frac{(1-aq^{2n})}{1-a}
\frac{\left(\frac{qa^2}{k}\right)_{2n}}{(k)_{2n}}\left(\frac{k}{a}\right)^n
\frac{\left(\rho_1,\rho_2,\frac{k}{\rho_1\rho_2}\right)_n} 
{\left(\frac{aq}{\rho_1},\frac{aq}{\rho_2},\frac{aq\rho_1\rho_2}{k}\right)_n}
\label{eq:2.10}
\end{equation}
and
\begin{equation}
\tilde\beta_n(a,k)=
\sum_{j=0}^n\frac{\left(\frac{k^2}{qa^2}\right)_{n-j}}{(q)_{n-j}}
\left(\frac{k^2}{qa^2}\right)^j 
\frac{\left(\frac{aq\rho_1}{k},\frac{aq\rho_2}{k},
\frac{aq}{\rho_1\rho_2},\frac{qa^2}{k}\right)_j}
{\left(\frac{aq}{\rho_1},\frac{aq}{\rho_2},\frac{aq\rho_1\rho_2}{k},q\right)_j}.
\label{eq:2.11}
\end{equation}
If we substitute this pair into \eqn{1.6}, the resulting identity is equivalent to 
Bailey nearly-poised transformation of a very-well-poised $_{12}\phi_{11}$ into a nearly-poised 
$_5\phi_4$ of the second kind ((III.25) in \cite{GR})
$$
_{12}W_{11}\left(a;a\sqrt{\frac{q}{k}},-a\sqrt{\frac{q}{k}},\frac{aq}{\sqrt k},-\frac{aq}{\sqrt k},
\rho_1,\rho_2,\frac{k}{\rho_1\rho_2},kq^n,q^{-n};q,q\right)=
$$
\begin{equation}
\frac{\left(aq,\frac{k^2}{qa^2}\right)_n}{\left(k,\frac{k}{a}\right)_n}
\mbox{ }_5\phi_4\left(\begin{array}{c} 
\frac{qa^2}{k},\frac{aq\rho_1}{k},\frac{aq\rho_2}{k},
\frac{aq}{\rho_1\rho_2},q^{-n};q,q \\
\frac{aq}{\rho_1},\frac{aq}{\rho_2},\frac{aq\rho_1\rho_2}{k},\frac{a^2q^{2-n}}{k^2}
\end{array}\right).
\label{eq:2.12}
\end{equation}

We conclude this section by remarking that if we continue a WP-tree further, 
we produce, in general, identities involving multisums. However, there is a way to generate 
further identities of the classical type involving only single fold sums.
Namely, we can fix free parameters in \eqn{2.9} and \eqn{2.12} in such a way that these 
transformations become summation formulas. In particular, if we replace $q,a,k$ by 
$\sqrt q,\sqrt a,\sqrt k$, respectively, and then set $\rho_1=i(q^2a)^{1/4},
\rho_2=-i(q^2a)^{1/4},\sigma_1=\sqrt kq^{n/2},\sigma_2=q^{-n/2}$ 
in \eqn{2.9}, then we get the summation formula:
$$
_{10}W_9\left(\sqrt a;i(q^2a)^{1/4},-i(q^2a)^{1/4},\frac{a}{k},
\sqrt kq^{n/2},-\sqrt kq^{n/2},q^{-n/2},-q^{-n/2};\sqrt q,\sqrt q\right)=
$$
\begin{equation}
\frac{\left(aq,\frac{a}{k};q\right)_n}{\left(\frac{k}{a},q\frac{k^2}{a};q\right)_n}
\frac{\left(-k\sqrt{\frac{q}{a}};\sqrt q\right)_{2n}}
{(-\sqrt a;\sqrt q)_{2n}}\left(\frac{k}{a\sqrt q}\right)^n.
\label{eq:2.13}
\end{equation}
This formula will play an important role in Section~4.

\medskip

\section{Further Bailey and Bailey-type transformations}
\label{sec:3}

To derive further $q$-hypergeometric transformations, 
we begin by setting $\rho_1=\frac{k}{aq}$ in formula \eqn{2.12}.
This yields
\begin{equation}
_{10}W_9\left(a;a\sqrt{\frac{q}{k}},-a\sqrt{\frac{q}{k}},\frac{aq}{\sqrt k},
-\frac{aq}{\sqrt k},\frac{k}{aq},kq^n,q^{-n};q,q\right)=
\frac{\left(aq,\frac{k^2}{qa^2}\right)_n}{\left(k,\frac{k}{a}\right)_n}.
\label{eq:3.1}
\end{equation}
In passing we note that the formula \eqn{3.1} is equivalent to the formula 
in Ex.~2.12 in \cite{GR}. Next, we observe that \eqn{3.1} contains a very-well-poised 
$_{10}W_9$ to which Bailey's transformation \eqn{2.9} applies. Consequently, after relabeling, 
we deduce from \eqn{3.1} that 
\begin{equation}
_{10}W_9\left(a;a\sqrt{\frac{q}{k}},-a\sqrt{\frac{q}{k}},\frac{a}{\sqrt k},
-\frac{aq}{\sqrt k},\frac{k}{a},kq^n,q^{-n};q,q\right)=
\frac{\left(aq,\sqrt k,\frac{k^2}{a^2}\right)_n}{\left(k,\frac{k}{a},q\sqrt k\right)_n}.
\label{eq:3.2}
\end{equation}
Identities \eqn{3.1} and \eqn{3.2} amount to the assertion that both $(\alpha^{(1)}_n,
\beta^{(1)}_n)$ and $(\alpha^{(2)}_n,\beta^{(2)}_n)$ form a WP-pair, where  

\begin{equation}
\alpha^{(1)}_n(a,k)=\frac{(a,q\sqrt a,-q\sqrt a)_n}{(q,\sqrt a,-\sqrt a)_n}
\frac{\left(\frac{qa^2}{k}\right)_{2n}}{(k)_{2n}}
\frac{\left(\frac{k}{aq}\right)_n}{\left(\frac{a^2q^2}{k}\right)_n}
\left(\frac{k}{a}\right)^n,
\label{eq:3.3}
\end{equation}

\begin{equation}
\beta^{(1)}_n(a,k)=\frac{\left(\frac{k^2}{qa^2}\right)_n}{(q)_n},
\label{eq:3.4}
\end{equation}
while
\begin{equation}
\alpha^{(2)}_n(a,k)=\frac{\left(a,q\sqrt a,-q\sqrt a,a\sqrt{\frac{q}{k}},-a\sqrt{\frac{q}{k}},
\frac{a}{\sqrt k},-\frac{aq}{\sqrt k},\frac{k}{a}\right)_n}
{\left(q,\sqrt a,-\sqrt a,\sqrt{qk},-\sqrt{qk},q\sqrt k,-\sqrt k,\frac{qa^2}{k}\right)_n},
\label{eq:3.5}
\end{equation}

\begin{equation}
\beta^{(2)}_n(a,k)=\frac{\left(\sqrt k,\frac{k^2}{a^2}\right)_n}{(q,q\sqrt k)_n}.
\label{eq:3.6}
\end{equation}
If we apply the construct \eqn{1.7}, \eqn{1.8} to $(\alpha^{(1)}_n,\beta^{(1)}_n)$ and 
$(\alpha^{(2)}_n,\beta^{(2)}_n)$ we find

$$
_{12}W_{11}\left(a;a\sqrt{\frac{q}{c}},-a\sqrt{\frac{q}{c}},\frac{aq}{\sqrt c},
-\frac{aq}{\sqrt c},\frac{c}{aq},\rho_1,\rho_2,kq^n,q^{-n};q,q\right)=
$$
\begin{equation}
\frac{\left(\frac{k\rho_1}{a},\frac{k\rho_2}{a},\frac{k}{c},aq\right)_n}
{\left(\frac{aq}{\rho_1},\frac{aq}{\rho_2},\frac{k}{a},cq\right)_n}
\mbox{}_7\phi_6\left(\begin{array}{c}
\frac{c^2}{qa^2},q\sqrt c,-q\sqrt c,\rho_1,\rho_2,kq^n,q^{-n};q,q \\
\sqrt c,-\sqrt c,\frac{k\rho_1}{a},\frac{k\rho_2}{a},\frac{cq^{1-n}}{k},cq^{1+n}
\end{array}\right)
\label{eq:3.7}
\end{equation}
and

$$
_{12}W_{11}\left(a;a\sqrt{\frac{q}{c}},-a\sqrt{\frac{q}{c}},\frac{a}{\sqrt c},
-\frac{aq}{\sqrt c},\frac{c}{a},\rho_1,\rho_2,kq^n,q^{-n};q,q\right)=
$$
\begin{equation}
\frac{\left(\frac{k\rho_1}{a},\frac{k\rho_2}{a},\frac{k}{c},aq\right)_n}
{\left(\frac{aq}{\rho_1},\frac{aq}{\rho_2},\frac{k}{a},cq\right)_n}
\mbox{}_6\phi_5\left(\begin{array}{c}
\frac{c^2}{a^2},-q\sqrt c,\rho_1,\rho_2,kq^n,q^{-n};q,q \\
-\sqrt c,\frac{k\rho_1}{a},\frac{k\rho_2}{a},\frac{c}{k}q^{1-n},cq^{1+n}
\end{array}\right),
\label{eq:3.8}
\end{equation}
where $c=\frac{k\rho_1\rho_2}{aq}$ as before. 
Employing formula \eqn{1.19}, one can show that \eqn{3.7} is equivalent to the famous 
$_{12}W_{11}\rightarrow\mbox{}_7\phi_6$ transformation of Bailey ((III.28) in \cite{GR}).
On the other hand, rewriting \eqn{3.8} with the aid of \eqn{1.19} results in the formula, 
which is equivalent to that of Jain (see Ex~2.14(ii), p. 52 in \cite{GR}).

To our inital surprise, if we insert $(\alpha^{(1)}_n,\beta^{(1)}_n)$ into \eqn{1.10} and 
\eqn{1.11}, the result is 

\begin{equation}
\alpha^{(3)}_n(a,k)=\alpha_n(a,k)
\label{eq:3.9}
\end{equation}
and

$$
\beta^{(3)}_n(a,k)=\sum_{j=0}^n
\frac{\left(\frac{k^2}{qa^2}\right)_{n-j}}{(q)_{n-j}}
\frac{\left(\frac{qa^2}{k^2}\right)_j}{(q)_j}
\left(\frac{k^2}{qa^2}\right)^j=
\frac{\left(\frac{k^2}{qa^2}\right)_n}{(q)_n}
\mbox{}_2\phi_1\left(\begin{array}{c}
q^{-n},\frac{qa^2}{k^2};q,q \\
\frac{a^2q^{2-n}}{k^2}\end{array}\right)
$$
\begin{equation}
\stackrel{\mbox{\scriptsize{by (II.6) in \cite{GR}}}}{=}
\frac{\left(\frac{k^2}{qa^2},q^{1-n}\right)_n}{\left(q,\frac{a^2q^{2-n}}{k^2}\right)_n}
\left(\frac{qa^2}{k^2}\right)^n=\delta_{n,0}=\beta_n(a,k).
\label{eq:3.10}
\end{equation}

\noindent
Thus we are back where we started in Section~2. Actually, in a view of a remark 
following \eqn{1.11}, this repetition of results should have been anticipated.

We now follow another lead suggested by Bailey (\cite{Ba1}, \S 10) wherein one sets half 
the $\alpha_n(a,k)$ to zero. In this way, we find the following WP--Bailey pair:

\begin{align}
\alpha_n^{(4)}(a,k)=
\begin{cases} 0, & \mbox{if } n \mbox{ is odd}, \\
              \left(\frac{k}{a}\right)^n
              \frac{\left(a,q^2\sqrt a,-q^2\sqrt a,\frac{a^2}{k^2};q^2\right)_{\frac{n}{2}}}
                   {\left(q^2,\sqrt a,-\sqrt a,\frac{q^2k^2}{a^2};q^2\right)_{\frac{n}{2}}},
                 & \mbox{if } n \mbox{ is even}, \end{cases}
\label{eq:3.11}
\end{align}
so to satisfy \eqn{1.6},

\begin{equation}
\beta^{(4)}_n(a,k)=\sum_{j\geq 0}\alpha_{2j}^{(4)}(a,k)
\frac{\left(\frac{k}{a}\right)_{n-2j}}{(q)_{n-2j}}
\frac{(k)_{n+2j}}{(aq)_{n+2j}}=
\frac{\left(k,\frac{k}{a}\right)_n}{(q,aq)_n}\times
\label{eq:3.12}
\end{equation}
$$
_8W_7\left(a;\frac{a^2}{k^2},kq^n,kq^{1+n},q^{-n},q^{1-n};q^2,q^2\right)=
\left(\frac{-k}{a}\right)^n \frac{\left(k,k\sqrt{\frac{q}{a}},-k\sqrt{\frac{q}{a}},
\frac{a}{k}\right)_n}
{\left(q,\sqrt{aq},-\sqrt{aq},q\frac{k^2}{a}\right)_n}
$$
by Jackson's $q$-analog of Dougall's theorem \eqn{2.5}.

We may now insert $(\alpha^{(4)}_n,\beta^{(4)}_n)$ pair into \eqn{1.7} and \eqn{1.8} 
to obtain a new WP-pair, and the result of putting this new pair into \eqn{1.6} is 
equivalent to the following Bailey-type identity

$$
_{12}W_{11}\left(a;\frac{a^2}{c^2},\rho_1,q\rho_1,\rho_2,q\rho_2,
kq^{1+n},kq^n,q^{1-n},q^{-n};q^2,q^2\right)=
$$
\begin{equation}
\frac{\left(aq,\frac{k}{c},\frac{k\rho_1}{a},\frac{k\rho_2}{a}\right)_n}
     {\left(cq,\frac{k}{a},\frac{aq}{\rho_1},\frac{aq}{\rho_2}\right)_n}
\mbox{}_{10}W_9\left(c;\rho_1,\rho_2,c\sqrt{\frac{q}{a}},-c\sqrt{\frac{q}{a}},
 \frac{a}{c},kq^n,q^{-n};q,\frac{-cq}{a}\right),
\label{eq:3.13}
\end{equation}
where, as before, $c=\frac{k\rho_1,\rho_2}{aq}$.

We note that \eqn{3.13} describes transformation of series with base $q^2$ to series 
with base $q$. While quite a few transformations of this kind can be found in the literature 
(see, for example, Sect.~3.10 in \cite{GR}), our result \eqn{3.13} appears to be new. 

Finally, we insert $(\alpha^{(4)}_n,\beta^{(4)}_n)$ into \eqn{1.10} and \eqn{1.11} 
to obtain another new WP-pair, and the result of putting this new pair into \eqn{1.6} is 
equivalent to the following formula

$$
\scriptstyle{
_{16}W_{15}\left(a;a\sqrt{\frac{q}{k}},-a\sqrt{\frac{q}{k}},\frac{aq}{\sqrt k},
-\frac{aq}{\sqrt k},\frac{aq^{\frac{3}{2}}}{\sqrt k},-\frac{aq^{\frac{3}{2}}}{\sqrt k},
\frac{aq^2}{\sqrt k},-\frac{aq^2}{\sqrt k},\frac{k^2}{a^2q^2},kq^{1+n},
kq^n,q^{1-n},q^{-n};q^2,q^2\right)}
$$
\begin{equation}
=\frac{\left(aq,\frac{k^2}{qa^2}\right)_n}{\left(k,\frac{k}{a}\right)_n}
\mbox{}_5\phi_4\left(\begin{array}{c}
\frac{a^2q}{k},\sqrt{\frac{a^3q^3}{k^2}},-\sqrt{\frac{a^3q^3}{k^2}},
\frac{k}{aq},q^{-n};q,\frac{-aq^2}{k} \\
\sqrt{aq},-\sqrt{aq},\frac{a^3q^3}{k^2},\frac{a^2q^{2-n}}{k^2}
\end{array}\right),
\label{eq:3.14}
\end{equation}
which also appears to be new. 

\medskip

\section{Bressoud's WP-pairs and their consequences}
\label{sec:4}

In \cite{Br}, Bressoud presented a striking variation on the classical work of Bailey and 
Rogers. Also in \cite{Br}, he found (what we may now call) three WP-Bailey pairs. 
This enabled him to discover new polynomial versions of the Rogers--Ramanujan identities:

\begin{equation}
\sum_{i=0}^Nq^{i^2}\qBin{N}{i}{q}=\sum_{j=-\infty}^\infty(-1)^jq^{\frac{j(5j+1)}{2}}
\qBin{2N}{N+2j}{q}  
\label{eq:4.1}
\end{equation}
and
 
\begin{equation}
\sum_{i=0}^{N}q^{i^2+i}\qBin{N}{i}{q}=\frac{1}{1-q^{N+1}}\sum_{j=-\infty}^\infty
(-1)^jq^{\frac{j(5j+3)}{2}}\qBin{2N+2}{N+2j+2}{q},  
\label{eq:4.2}
\end{equation}
where

\begin{align}
\qBin{N}{j}{q}
=\begin{cases} \frac{(q^{N-j+1})_j}{(q)_j}, & \mbox{if } j\geq 0, \\
               0, & \mbox{otherwise,} \end{cases}
\label{eq:4.3}
\end{align}
along with finite versions of the Rogers--Selberg identities.

Bressoud's variation of Bailey Lemma can be stated as:

\begin{thm}\textup{(}Bressoud\textup{)}
If $(\alpha_n(a,k),\beta_n(a,k))$ form WP-Bailey pair, then the following identity holds true

\begin{equation}
\sum_{j\geq 0}\beta_j\left(a,q\frac{a^2}{k}\right)\left(\frac{k^2}{qa^2}\right)^j=
\frac{\left(k,\frac{k}{a}\right)_\infty}{\left(aq,\frac{k^2}{qa^2}\right)_\infty}
\sum_{j\geq 0}\frac{\left(\frac{qa^2}{k}\right)_{2j}}{(k)_{2j}}
\left(\frac{k^2}{qa^2}\right)^j\alpha_j\left(a,q\frac{a^2}{k}\right).
\label{eq:4.4}
\end{equation}

\label{thm:1}
\end{thm}
It is easy to recognize this theorem as the limiting case of our construct \eqn{1.10} and 
\eqn{1.11} with $n\rightarrow\infty$ in \eqn{1.11}. In addition, Bressoud's first WP-pair 
is equivalent to our $(\alpha'_n,\beta'_n)$ of Section~2 with $\rho_2\rightarrow\infty$. 
These two observations imply that equation \eqn{3.2} of \cite{Br} and, as a result, 
Bressoud's identities \eqn{4.1} and \eqn{4.2} are, actually, the limiting special cases of
Bailey's $_{12}W_{11}\rightarrow\mbox{}_5\phi_4$ transformation formula \eqn{2.12}. 
Surprisingly, this $q$-hypergeometric explanation of Bressoud's formulas has never been 
given before.

Bressoud's second WP-Bailey pair is given by 

\begin{equation}
\alpha_n^*(a,k)=\frac{1-\sqrt aq^n}{1-\sqrt a}
\frac{\left(\sqrt a,\sqrt q\frac{a}{k};\sqrt q\right)_n}
     {\left(\sqrt q,\frac{k}{\sqrt a};\sqrt q\right)_n}
\left(\frac{k}{a\sqrt q}\right)^n
\label{eq:4.5}
\end{equation}
and

\begin{equation}
\beta_n^*(a,k)=\frac{\left(k,\frac{aq}{k}\right)_n}{\left(q,\frac{k^2}{a}\right)_n}
\frac{\left(-\frac{k}{\sqrt a};\sqrt q\right)_{2n}}
     {\left(-\sqrt{aq};\sqrt q\right)_{2n}}
\left(\frac{k}{a\sqrt q}\right)^n.
\label{eq:4.6}
\end{equation}
The assertion that \eqn{1.6} holds for this pair is  special case of Jackson's 
$q$-analog of Dougall's theorem \eqn{2.5}.

Bressoud's third and last WP-pair is

\begin{equation}
\alpha_n^\dagger(a,k)=\frac{1-\sqrt aq^{2n}}{1-a}
\frac{\left(\sqrt a,\frac{a}{k};\sqrt q\right)_n}
     {\left(\sqrt q,k\sqrt{\frac{q}{a}};\sqrt q\right)_n}
\left(\frac{k}{a\sqrt q}\right)^n
\label{eq:4.7}
\end{equation}
and

\begin{equation}
\beta_n^\dagger(a,k)=\frac{\left(k,\frac{a}{k},-k\sqrt{\frac{q}{a}},-\frac{kq}{\sqrt a}\right)_n}
{\left(q,\frac{qk^2}{a},-\sqrt a,-\sqrt{aq}\right)_n}
\left(\frac{k}{a\sqrt q}\right)^n.
\label{eq:4.8}
\end{equation}
Bressoud's proof that \eqn{4.7}, \eqn{4.8} satisfy \eqn{1.6} is the most difficult of 
the three pairs he considers. It combines a double series expansion of a very-well-poised 
$_{10}W_9$ along with an application of the $q$--Pfaaf\--Saalsch\"utz summation ((II.12) in
\cite{GR}). However, it turns out that his argument may be simplified greatly by the 
observation that his final $_{10}W_9$ summation formula is nothing else but \eqn{2.13}: 
a special case of Bailey's $_{10}W_9\rightarrow\mbox{}_{10}W_9$ transformation. In \cite{Br}, 
pairs $(\alpha_n^*,\beta_n^*)$ and $(\alpha_n^\dagger,\beta_n^\dagger)$ were employed to derive polynomial 
versions of the Rogers--Selberg identities.

Let us now move on to discuss $q$-hypergeometric consequences of Bressoud's second and third 
pairs. To this end we create new WP-pairs, using  $(\alpha_n^*,\beta_n^*)$ and 
$(\alpha_n^\dagger,\beta_n^\dagger)$ as input, via \eqn{1.7} and \eqn{1.8} or \eqn{1.10} and 
\eqn{1.11}. The four identities below are the result of putting these new pairs into \eqn{1.6}. 
In order to improve the appearance of these identities we replace $a,q,k,\rho_1,\rho_2$ by
$a^2,q^2,\rho_1^2,\rho_2^2$, respectively.\\
From the Bressoud's second pair $(\alpha_n^*,\beta_n^*)$ we deduce

\begin{equation}
_{10}W_9\left(c^2;\rho_1^2,\rho_2^2,\frac{a^2q^2}{c^2},-\frac{c^2}{a},-\frac{qc^2}{a},
k^2q^{2n},q^{-2n};q^2,q\frac{c^2}{a^2}\right)=
\label{eq:4.9}
\end{equation}
$$
\frac{\left(\frac{a^2q^2}{\rho_1^2},\frac{a^2q^2}{\rho_2^2},c^2q^2,\frac{k^2}{a^2};q^2\right)_n}
     {\left(\frac{k^2\rho_1^2}{a^2},\frac{k^2\rho_2^2}{a^2},a^2q^2,\frac{k^2}{c^2};q^2\right)_n}
\mbox{}_{12}W_{11}
\left(a;\rho_1,-\rho_1,\rho_2,-\rho_2,\frac{qa^2}{c^2},kq^n,-kq^n,q^{-n},-q^{-n};q,q\right),
$$
where $c$ is defined in \eqn{1.9} and

$$
\Phi\begin{bmatrix}
a,q\sqrt a,-q\sqrt a,\frac{k^2}{qa^2}:q^{-2n},k^2q^{2n},\frac{qa^2}{k},-\frac{qa^2}{k},
\frac{q^2a^2}{k},-\frac{q^2a^2}{k};q,q^2;q \\
\sqrt a,-\sqrt a,\frac{a^3q^2}{k^2}:\frac{a^2q^{2-2n}}{k^2},a^2q^{2+2n},k,-k,kq,-kq
\end{bmatrix}=
$$
\begin{equation}
\frac{\left(a^2q^2,\frac{k^4}{q^2a^4};q^2\right)_n}
     {\left(k^2,\frac{k^2}{a^2};q^2\right)_n}
\mbox{}_5\phi_4\begin{pmatrix}
\frac{a^4q^2}{k^2},\frac{k^2}{a^2},-\frac{a^3q^2}{k^2},-\frac{a^3q^3}{k^2},q^{-2n};
      q^2,\frac{a^2q^3}{k^2} \\
\frac{a^6q^4}{k^4},-aq,-aq^2,\frac{a^4q^{4-2n}}{k^4}
\end{pmatrix}.
\label{eq:4.10}
\end{equation}
From Bressoud's third pair $(\alpha_n^\dagger,\beta_n^\dagger)$ we deduce

\begin{equation}
_{10}W_9\left(c^2;\rho_1^2,\rho_2^2,\frac{a^2}{c^2},-c^2\frac{q}{a},-\frac{c^2q^2}{a},
k^2q^{2n},q^{-2n};q^2,\frac{c^2}{a^2}q\right)=
\label{eq:4.11}
\end{equation}
$$
\frac{\left(\frac{a^2q^2}{\rho_1^2},\frac{a^2q^2}{\rho_2^2},c^2q^2,\frac{k^2}{a^2};q^2\right)_n}
     {\left(\frac{k^2\rho_1^2}{a^2},\frac{k^2\rho_2^2}{a^2},a^2q^2,\frac{k^2}{c^2};q^2\right)_n}
\Phi\begin{bmatrix}
a,\frac{a^2}{c^2}:\rho_1^2,\rho_2^2,aq^2,-aq^2,k^2q^{2n},q^{-2n};q,q^2;q \\
c^2\frac{q}{a}:\left(\frac{aq}{\rho_1}\right)^2,\left(\frac{aq}{\rho_2}\right)^2,a,-a,
\left(\frac{a}{k}\right)^2q^{2-2n},a^2q^{2+2n}
\end{bmatrix},
$$
with $c$ as in \eqn{1.9} and

$$
\Phi\begin{bmatrix}
a,\frac{k^2}{a^2q^2}:q^2a,-q^2a,\frac{a^2q}{k},-\frac{a^2q}{k},
\frac{a^2q^2}{k},-\frac{a^2q^2}{k},k^2q^{2n},q^{-2n};q,q^2;q \\
\frac{a^3q^3}{k^2}:a,-a,k,-k,qk,-qk,\frac{a^2q^{2-2n}}{k^2},a^2q^{2+2n}
\end{bmatrix}=
$$
\begin{equation}
\frac{\left(\frac{k^4}{a^4q^2},a^2q^2;q^2\right)_n}
     {\left(k^2,\frac{k^2}{a^2};q^2\right)_n}
\mbox{}_5\phi_4\begin{pmatrix}
\frac{a^4q^2}{k^2},\frac{k^2}{a^2q^2},-\frac{a^3q^3}{k^2},-\frac{a^3q^4}{k^2},q^{-2n};
      q^2,\frac{a^2q^3}{k^2} \\
\frac{a^6q^6}{k^4},-a,-aq,\frac{q^{4-2n}}{k^4}a^4\end{pmatrix}.
\label{eq:4.12}
\end{equation}
Formulas \eqn{4.9}-\eqn{4.12} appear to be new.

\medskip

\section{Doubly bounded polynomial identities of the Rogers--Ramanujan type}
\label{sec:5}

As we mentioned at the beginning of Section~4, Bressoud's noteworthy achievement in 
\cite{Br} was a discovery and proof of two new polynomial identities \eqn{4.1} and \eqn{4.2},
which in the limit as $N\rightarrow\infty$ yield Rogers--Ramanujan identities.

Now that we have applied his WP--Bailey pairs to produce the successor WP-pairs in 
the WP--Bailey tree, we may deduce a number of doubly bounded polynomial identities 
of the Rogers--Ramanujan type.

There is ample precedent for such results \cite{Bu}, \cite{F}, \cite{War2}. 
In particular, Burge \cite{Bu} presented a doubly bounded identity that 
contained both the finite version of the Rogers--Ramanujan identity given in \cite{An4} as 
well as the Schur's celebrated polynomial analog of the Rogers--Ramanujan identity \cite{Sch}. 
Recently, Warnaar \cite{War2}, drawing on work of Burge \cite{Bu} and Foda et al \cite{F}, 
also found two parameter polynomial identities, which not only contain results of \cite{An4}, 
but also \eqn{4.1} as limiting cases:

\begin{equation}
\sum_{i=0}^M q^{i^2}\qBin{N}{i}{q}\qBin{2N+M-i}{2N}{q}=
\sum_{j=-M}^M (-1)^jq^{\frac{j(5j-1)}{2}}
\qBin{N+M+j}{N+2j}{q}\qBin{N+M-j}{N-2j}{q}.  
\label{eq:5.1}
\end{equation}
In the same work, Warnaar finds a similar formula for the second 
Rogers--Ramanujan identity:

\begin{equation}
\sum_{i=0}^{M-1}q^{i^2+i}\qBin{N}{i}{q}\qBin{2N+M-i}{2N+1}{q}=
\sum_{j=-M}^{M-1}(-1)^jq^{\frac{j(5j+3)}{2}}
\qBin{N+M+j}{N+2j+1}{q}\qBin{N+M-j}{N-2j}{q},  
\label{eq:5.2}
\end{equation}
where here and throughout $N,M\in \Bbb{Z}_{\geq 0}$. 
We note that \eqn{5.1} follows from Bailey's $_{12}W_{11}\rightarrow\mbox{}_5\phi_4$ 
transformation formula \eqn{2.12} under the substitutions $n=M,a=1,k=q^{1+N},
\rho_1\rightarrow\infty,\rho_2\rightarrow\infty$. If we set $n=M,a=q,k=q^{3+N},
\rho_1\rightarrow\infty,\rho_2\rightarrow\infty$ in \eqn{2.12}, we get a formula, which 
is very similar to \eqn{5.2}:

\begin{equation}
\sum_{i=0}^M q^{i^2+i}\qBin{N}{i}{q}\qBin{2N+M-i+2}{2N+2}{q}=
\label{eq:5.3}
\end{equation}
$$
\frac{1}{1-q^{N+1}}
\sum_{j=-M-1}^M(-1)^jq^{\frac{j(5j+3)}{2}}
\qBin{N+M+j+2}{N+2j+2}{q}\qBin{N+M-j+1}{N-2j}{q}.  
$$
We observe that \eqn{5.3} contains Bressoud's formula \eqn{4.2} as a limiting case 
$M\rightarrow\infty$.

Just as there are numerous important classical identities contained in Watson's 
$q$-analog of Whipple's theorem (eq.(III.17) in \cite{GR}), so too there are further 
equally interesting special cases of \eqn{2.12}.
In particular, setting $n=M,a=1,k=q^{1+N},
\rho_1=-1,\rho_2\rightarrow\infty$ in \eqn{2.12} yields

\begin{equation}
\sum_{i=0}^M q^{i^2}\qBin{N}{i}{q^2}\qBin{2N+M-i}{2N}{q}=
\sum_{j=-M}^M(-1)^jq^{2j^2}
\qBin{N+M-j}{N-2j}{q}\qBin{N+M+j}{N+2j}{q},  
\label{eq:5.4}
\end{equation}
a two parameter polynomial generalization of

\begin{equation}
\sum_{i\geq 0}\frac{q^{i^2}}{(q^2,q^2)_i}=\frac{1}{(q)_\infty}
\sum_{j=-\infty}^\infty(-1)^jq^{2j^2}.  
\label{eq:5.5}
\end{equation}
We note that \eqn{5.4} is the $(a,b)=(2,1)$ case of Theorem 7.1 in \cite{War2}.

Now, set $n=M,a=1,k=q^{1+N},\rho_1=x,\rho_2=\frac{1}{x^2}$ in \eqn{2.12} and let 
$x\rightarrow 0$ to find another formula in \cite{War2}

\begin{equation}
\sum_{i=0}^M q^{Ni}\qBin{N}{i}{q}\qBin{2N+M-i}{2N}{q}=
\sum_{j=-M}^M(-1)^jq^{\frac{3j^2-j}{2}}
\qBin{N+M+j}{N+2j}{q}\qBin{N+M-j}{N-2j}{q}.  
\label{eq:5.7}
\end{equation}
This is a two parameter refinement of Bressoud's finite version of Euler's Pentagonal Number 
Theorem \cite{Br}. Not surprisingly, identity \eqn{5.7} reduces to \eqn{5.1} under the 
substitution $q\rightarrow q^{-1}$.

Next, replace $q$ by $q^2$ in \eqn{2.12}. Then, set $n=M,a=1,k=q^{2+2N},
\rho_1=q,\rho_2\rightarrow\infty$. This yields a two parameter refinement of the 
Andrews--Santos identity \cite{AS}:

\begin{equation}
\sum_{i=0}^M q^{2i^2}\qBin{2N}{2i}{q}\qBin{2N+M-i}{2N}{q^2}=
\sum_{j=-M}^M q^{4j^2-j}
\qBin{N+M+j}{N+2j}{q^2}\qBin{N+M-j}{N-2j}{q^2}.  
\label{eq:5.6}
\end{equation}
In the limit as $N,M\rightarrow\infty$, the above identity reduces to eqn. (39) on 
the celebrated Slater's list \cite{S}. 

If we set $n=M,a=1,k=q^{1+N},\rho_1=\rho_2=-1$ in \eqn{2.12} we obtain

\begin{equation}
\sum_{i=0}^{\min (N,M)} q^{i^2}\qBin{2N+M-i}{2N}{q}
\frac{\qBinpwr{N}{i}{q^2}{2}}{\qBinpwr{N}{i}{q}{2}}=
\label{eq:5.8}
\end{equation}
$$
2\sum_{j=-\min(M,\lfloor\frac{N}{2}\rfloor)}^{\min(M,\lfloor\frac{N}{2}\rfloor)}
(-1)^j\frac{q^{\frac{3j+1}{2}j}}{1+q^j}
\qBin{N+M+j}{N}{q}\qBin{N+M-j}{N}{q}\frac{\qBin{2N}{N+2j}{q}}{\qBin{2N}{N+j}{q}}.  
$$
In the limit as $N\rightarrow\infty,M\rightarrow\infty$, this is revealed to be 
Watson's formula \cite{Wat} for the third order mock theta function:

\begin{equation}
\sum_{i\geq 0} \frac{q^{i^2}}{(-q)_i^2}=\frac{2}{(q)_\infty}
\sum_{j=-\infty}^\infty(-1)^j\frac{q^{\frac{3j+1}{2}j}}{1+q^j}.
\label{eq:5.9}
\end{equation}
Or one may take $n=M,a=1,k=q^{1+N},\rho_1=1,\rho_2=\infty$ in \eqn{2.12} to obtain 

\begin{equation}
\sum_{i=0}^M q^{i^2}\qBinpwr{N}{i}{q}{2}\qBin{2N+M-i}{2N}{q}=
\qBinpwr{N+M}{N}{q}{2},
\label{eq:5.10}
\end{equation}
which is a two parameter generalization of the classic expansion of the generating function 
for the unrestricted partitions (eqn.~(2.2.9) in \cite{An5}):

\begin{equation}
\sum_{i\geq 0}\frac{q^{i^2}}{(q)_i^2}=\frac{1}{(q)_\infty}.  
\label{eq:5.11}
\end{equation}
In fact, \eqn{5.10} is a special case of the $q$-Pfaff-Saalsch\"utz identity 
((II.12) in \cite{GR}).

While the formulas of this section have so far been derived from Bailey's 
$_{12}W_{11}\rightarrow\mbox{}_5\phi_4$ transformation formula \eqn{2.12}, it is definitely
possible to obtain further doubly bounded identities from some of the other transformation 
formulas. For example, if we set $n=M,a=1,k=q^{1+N}$ in \eqn{4.10}, we obtain

\begin{equation}
\sum_{i=0}^{\min(N,M)} q^{2i^2}\qBin{N}{i}{q^2}\frac{\qBin{2N}{2i}{q^2}}{\qBin{2N}{2i}{q}}
\frac{\qBin{N+i}{N}{q^2}}{\qBin{2N}{i}{q^2}}\qBin{2N+M-i}{2N}{q^2}=
\frac{\qBin{4N}{2N}{q}}{\qBin{2N}{N}{q^2}}
\label{eq:5.12}
\end{equation}
$$
\times\sum_{j=-\min(M,\lfloor\frac{N}{2}\rfloor)}^{\min(M,\lfloor\frac{N}{2}\rfloor)}
(-1)^jq^{\frac{j(7j-1)}{2}}
\qBin{N+M+j}{N}{q^2}\qBin{N+M-j}{N}{q^2}\frac{\qBin{2N}{N-2j}{q^2}}{\qBin{4N}{2N-j}{q}},  
$$
and if we let $M$ and $N\rightarrow\infty$ we end up with one of the Rogers--Selberg identities 
(Ex.~10,p.117 in \cite{An5})

\begin{equation}
\sum_{i\geq 0}\frac{q^{2i^2}}{(q^2;q^2)_i(-q)_{2i}}=\frac{1}{(q^2,q^2)_\infty}
\sum_{j=-\infty}^\infty(-1)^jq^{\frac{j(7j-1)}{2}}.  
\label{eq:5.13}
\end{equation}
We note that \eqn{5.12} contains Bressoud's polynomial version of \eqn{5.13} as a special 
limiting case $M\rightarrow\infty$.

\medskip

\section{The Gasper--Kirillov identity}
\label{sec:6}

We continue this account of WP--Bailey pairs and their applications with a consideration of 
an identity communicated to one of us (G.E.A.) by A.N.~Kirillov. The attempt to prove Kirillov's 
identity from the existing literature led to the constructs \eqn{1.7}--\eqn{1.11}
\cite{An2}.
In praise of nescience we add that the entire effort would never have been made if we had 
observed that Kirillov's formula is, in fact, the case $\delta=-1$ of the following formula

$$
_{12}W_{11}\left(
a;a\sqrt{\frac{q}{k}},-a\sqrt{\frac{q}{k}},\frac{aq}{\sqrt k},-\frac{aq}{\sqrt k},
\frac{k}{a},\frac{k}{a\delta},a\delta,kq^n,q^{-n};q,\frac{aq}{k}\right)=
$$
\begin{equation}
\frac{\left(aq,\frac{k}{a\delta}\right)_n}
     {\left(k,\frac{q}{\delta}\right)_n}
\mbox{}_4\phi_3\begin{pmatrix}
\frac{k}{a},a\delta,q^{-n},\delta q^{-n};q,\left(\frac{aq}{k}\right)^2 \\
\delta\frac{a^2q}{k},\frac{aq^{1-n}}{k},\delta\frac{aq^{1-n}}{k}
\end{pmatrix},
\label{eq:6.1}
\end{equation}
which is equivalent to the formula of Gasper \cite{G},[\cite{GR}, p.~231, Ex.~8.15].

To prove \eqn{6.1} we use the following WP-pair:

\begin{equation}
\alpha_n^W(a,k)=
\frac{\left(a,\delta a,q\sqrt a,-q\sqrt a,\frac{aq}{k},\frac{aq}{k\delta}\right)_n}
     {\left(q,\frac{q}{\delta},\sqrt a,-\sqrt a,k,k\delta\right)_n}
\left(\frac{k^2}{qa^2}\right)^n
\label{eq:6.2}
\end{equation}
and

\begin{equation}
\beta_n^W(a,k)=\frac{\left(\frac{k}{a},\frac{k}{a\delta}\right)_n}
{\left(q,\frac{q}{\delta}\right)_n}
\mbox{}_4\phi_3\begin{pmatrix}
\delta q^{-n},q^{-n},\frac{aq}{k},a\delta;q,q \\
k\delta,\frac{a}{k}q^{1-n},\delta\frac{a}{k}q^{1-n}
\end{pmatrix}.
\label{eq:6.3}
\end{equation}
The assertion that $\alpha_n^W$ and $\beta_n^W$ satisfy \eqn{1.6} is an instance of 
Watson's $q$-analog of Whipple's theorem (eq.(III.17) in \cite{GR}).

If we now put this WP--Bailey pair into \eqn{1.10} and \eqn{1.11} we obtain a new 
WP-pair, which when inserted into \eqn{1.6} yields the following result:

$$
_{12}W_{11}\left(
a;a\sqrt{\frac{q}{k}},-a\sqrt{\frac{q}{k}},\frac{aq}{\sqrt k},-\frac{aq}{\sqrt k},
\frac{k}{a},\frac{k}{\delta a},\delta a,kq^n,q^{-n};q,\frac{aq}{k}\right)=
$$
\begin{equation}
\frac{\left(aq,\frac{k^2}{qa^2}\right)_n}
     {\left(k,\frac{k}{a}\right)_n}
\sum_{j=0}^n\frac{\left(q^{-n},\frac{aq}{k},\frac{aq}{k\delta}\right)_j}
                 {\left(q,\frac{q}{\delta},\frac{a^2q^{2-n}}{k^2}\right)_j}q^j
\mbox{}_4\phi_3\begin{pmatrix}
q^{-j},\delta q^{-j},\frac{k}{a},a\delta;q,q \\
\delta\frac{a^2q}{k},\frac{kq^{-j}}{a},\delta\frac{kq^{-j}}{a}
\end{pmatrix}.
\label{eq:6.4}
\end{equation}
Clearly \eqn{6.4} has the same left-hand side as \eqn{6.1}, but definitely not the same 
right-hand side. 
To identify these right-hand sides, we use the $q$-Pfaff-Saalsch\"utz sum ((II.12) in \cite{GR}) 
and note that the right-hand side of \eqn{6.4} may be written as
 
$$
\frac{\left(\frac{k^2}{qa^2},aq\right)_n}
     {\left(k,\frac{k}{a}\right)_n}
\sum_{i\geq 0}\frac{\left(\frac{k}{a},\delta a\right)_i}
                   {\left(q,\delta\frac{qa^2}{k}\right)_i}q^i
\sum_{j=0}^{n-i}q^{i+j}\frac{\left(q^{-i-j},\delta q^{-i-j}\right)_i}
                            {\left(\frac{k}{a}q^{-i-j},\delta\frac{k}{a}q^{-i-j}\right)_i}
\frac{\left(q^{-n},\frac{aq}{k},\frac{aq}{k\delta}\right)_{i+j}}
     {\left(q,\frac{q}{\delta},\frac{a^2q^{2-n}}{k^2}\right)_{i+j}}
$$
$$
=\frac{\left(\frac{k^2}{qa^2},aq\right)_n}
     {\left(k,\frac{k}{a}\right)_n}
\sum_{i=0}^n\frac{\left(\delta a,\frac{k}{a},q^{-n}\right)_i}
                 {\left(q,\delta\frac{qa^2}{k},\frac{a^2q^{2-n}}{k^2}\right)_i}
\left(\frac{aq}{k}\right)^{2i}
\mbox{}_3\phi_2\begin{pmatrix}
\frac{aq}{k},\frac{aq}{k\delta},q^{-n+i};q,q \\
\frac{q}{\delta},\frac{a^2q^{2-n+i}}{k^2}
\end{pmatrix}
$$
$$
\stackrel{\mbox{\scriptsize{by (II.12) in \cite{GR}}}}{=}
\frac{\left(\frac{k^2}{qa^2},aq\right)_n}
     {\left(k,\frac{k}{a}\right)_n}
\sum_{i=0}^n\frac{\left(\delta a,\frac{k}{a},q^{-n}\right)_i}
                 {\left(q,\delta\frac{qa^2}{k},\frac{a^2q^{2-n}}{k^2}\right)_i}
\left(\frac{aq}{k}\right)^{2i}
\frac{\left(\frac{k}{a},\frac{k}{a\delta}\right)_{n-i}}
     {\left(\frac{q}{\delta},\frac{k^2}{qa^2}\right)_{n-i}}
$$
\begin{equation}
=\frac{\left(aq,\frac{k}{a\delta}\right)_n}
      {\left(k,\frac{q}{\delta}\right)_n}
\mbox{}_4\phi_3\begin{pmatrix}
\frac{k}{a},\delta a,q^{-n},\delta q^{-n};q,\left(\frac{aq}{k}\right)^2 \\
\delta\frac{qa^2}{k},\frac{aq^{1-n}}{k},\delta\frac{aq^{1-n}}{k}
\end{pmatrix},
\label{eq:6.5}
\end{equation}
as desired.

\medskip

\section{Well-poised Burge pairs and further doubly bounded polynomial identities}
\label{sec:7}

We say that sequences $\{A_N\}$ and $\{B_N\}$ form a Burge pair if they satisfy the 
following relation

\begin{equation}
B_N(M,q)=\sum_{j\geq 0}A_j(q)Q(N,M,aj,bj,q)
\label{eq:7.1}
\end{equation}
with $N,M,a,b\in \Bbb{Z}_{\geq 0}$, and

\begin{equation}
Q(N,M,aj,bj,q)=\qBin{N+M+aj-bj}{N+aj}{q}\qBin{N+M-aj+bj}{N-aj}{q}.
\label{eq:7.2}
\end{equation}
In \cite{Bu}, Burge made a crucial use of the following formulas

\begin{equation}
q^{b^2j^2}Q(N,M,(a+b)j,bj,q)=
\sum_{i=0}^{M}q^{i^2}\qBin{2N+M-i}{2N}{q}Q(N-i,i,aj,bj,q),
\label{eq:7.3}
\end{equation}
and

\begin{equation}
q^{a^2j^2}Q(N,M,(a+b)j,aj,q)=\sum_{i=0}^{M}q^{i^2}\qBin{2N+M-i}{2N}{q}
Q(i,N-i,aj,bj,q),
\label{eq:7.4}
\end{equation}
with $N,M,a,b,\in\Bbb{Z}_{\geq 0}$ and $N\ge |(a-b)j|$. It is easy to check that 
\eqn{7.3} and \eqn{7.4} are immediate consequences of the $q$-Pfaff-Saalsch\"utz 
sum((II.12) in \cite{GR}).
The power of \eqn{7.3} and \eqn{7.4} lies in the fact that these transformations can 
be employed to generate an infinite binary (Burge) tree from the initial ``seed" 
identity \eqn{7.1}.
In particular, if, following \cite{Bu}, we apply \eqn{7.3} to \eqn{7.1} with $0\le b\le 2a$, 
we obtain the new identity

\begin{equation}
\sum_{i=0}^{M}q^{i^2}\qBin{2N+M-i}{2N}{q}B_{N-i}(i,q)=
\sum_{j\geq 0}A_j(q)q^{b^2j^2}Q(N,M,(a+b)j,bj,q).
\label{eq:7.5}
\end{equation}
Similar application of \eqn{7.4} to \eqn{7.1} with $0\le b\le 2a$ yields

\begin{equation}
\sum_{i=0}^{M}q^{i^2}\qBin{2N+M-i}{2N}{q}B_i(N-i,q)=
\sum_{j\ge 0} A_j(q)q^{a^2j^2}Q(N,M,(a+b)j,aj,q).
\label{eq:7.6}
\end{equation}
We now observe that both \eqn{7.5} and \eqn{7.6} are of the form \eqn{7.1}. Therefore, we can 
transform these identities into the four new ones and so it goes.

Comparison of \eqn{7.1} with $a=1,b=0$ and \eqn{1.6} with $n=N,k=q^{1+M},a=1$ suggests 
that a WP--Bailey pair $(\alpha_n(1,k), \beta_n(1,k))$ can be interpreted as a Burge pair, 
provided that $\alpha_n(1,k)$ does not depend on $k$. We will call such pairs WP--Burge pairs.

If we now insert \eqn{2.3} and \eqn{2.4} into \eqn{1.6} with $n=N$, and then set 
$a=1,k=q^{1+M},\rho_1=-1,\rho_2\rightarrow\infty$, we find

\begin{equation}
(-1)^N \qBin{N+M}{N}{q^2}=
\sum_{j=-\infty}^\infty(-1)^j Q(N,M,j,0,q).
\label{eq:7.7}
\end{equation}
Actually, this identity is a special case of the $q$-Kummer sum ((II.9) in \cite{GR}).
From \eqn{7.7} we read off our first WP--Burge pair:

\begin{align}
A_N^{(1)}(q)=
\begin{cases} 2(-1)^N, & \mbox{if } N>0, \\
              1,       & \mbox{if } N=0, \end{cases}
\label{eq:7.8}
\end{align}

\begin{equation}
B_N^{(1)}(M,q)=(-1)^N \qBin{N+M}{N}{q^2}.
\label{eq:7.9}
\end{equation}
Iterating \eqn{7.7} by \eqn{7.4} yields

\begin{equation}
\qBin{2N+M}{N}{q}
\mbox{}_3\phi_2\begin{pmatrix}
q^{-N},-q^{-N},q^{-M};q,q \\
-q,q^{-2N-M}
\end{pmatrix}=
\sum_{j=-M}^{M} (-1)^j q^{j^2} 
Q(N,M,j,j,q).
\label{eq:7.10}
\end{equation}
By evaluating the left-hand side with the aid of the $q$--Pfaff-Saalsch\"utz sum 
((II.12) in \cite{GR}), we arrive at Burge's identity \cite{Bu}

\begin{equation}
\qBin{N+M}{N}{q^2}=\sum_{j=-M}^{M} (-1)^j q^{j^2} 
Q(N,M,j,j,q).
\label{eq:7.11}
\end{equation}
Further iterations of \eqn{7.11} by Burge transforms \eqn{7.3} and \eqn{7.4} can be 
found in \cite{War2}.

Next, in \eqn{1.6} with \eqn{2.3}, \eqn{2.4} and $n=N$, we set $a=1,k=q^{1+M},
\rho_1\rightarrow\infty,\rho_2\rightarrow\infty$ to get 

\begin{equation}
(-1)^Nq^{\binom{N+1}{2}+NM}\qBin{N+M}{N}{q}=
\sum_{j=-\infty}^\infty (-1)^j q^{\binom{j}{2}} 
Q(N,M,j,0,q),
\label{eq:7.12}
\end{equation}
which is defining relation for our second WP--Burge pair:

\begin{align}
A_N^{(2)}(q)=
\begin{cases} (-1)^N(1+q^N)q^{\binom{N}{2}}, & \mbox{if } N>0, \\
              1,                             & \mbox{if } N=0, \end{cases}
\label{eq:7.13}
\end{align}

\begin{equation}
B_N^{(2)}(M,q)=(-1)^N \qBin{N+M}{N}{q}q^{\binom{N+1}{2}+NM}.
\label{eq:7.14}
\end{equation}
It is worth mentioning that \eqn{7.12} is, in fact, a limiting case of the Jackson's 
terminating $q$-analog of Dixon's sum ((II.15) in \cite{GR}).

If we apply \eqn{7.4} to \eqn{7.12} and use $q$-Chu-Vandermonde sum ((II.6) in \cite{GR}), 
then we rediscover a recent result of Warnaar \cite{War2}:

\begin{equation}
\qBin{N+M}{N}{q}=\sum_{j=-M}^{M} (-1)^j q^{\frac{3j+1}{2}j} 
Q(N,M,j,j,q).
\label{eq:7.15}
\end{equation}
Further iterations of \eqn{7.15} by \eqn{7.3} and \eqn{7.4} are described in 
detail in \cite{War2}. Here, we confine 
ourselves to the comment that the doubly bounded polynomial analog of the first 
Rogers--Ramanujan identity \eqn{5.1} is an immediate consequence of \eqn{7.15} and \eqn{7.6}.

To derive our third WP--Burge pair, we again use \eqn{1.6} with \eqn{2.3} and \eqn{2.4}.
This time we set $n=N,a=1,k=q^{1+M},\rho_1=\sqrt q,\rho_2\rightarrow\infty$ to obtain

\begin{equation}
q^{-\frac{N}{2}}\qBin{2N+2M+1}{2N}{\sqrt q}=
\sum_{j=-\infty}^\infty q^{\frac{j}{2}} Q(N,M,j,0,q).
\label{eq:7.16}
\end{equation}
Clearly, this identity amounts to the assertion that $A_N^{(3)}(q), B_N^{(3)}(M,q)$ 
form a WP--Burge pair, where

\begin{align}
A_N^{(3)}(q)=
\begin{cases} q^{\frac{N}{2}}+q^{-\frac{N}{2}}, & \mbox{if } N>0, \\
              1,                                & \mbox{if } N=0, \end{cases}
\label{eq:7.17}
\end{align}
and
\begin{equation}
B_N^{(3)}(M,q)=q^{-\frac{N}{2}} \qBin{2N+2M+1}{2N}{\sqrt q}.
\label{eq:7.18}
\end{equation}
To get a better feeling for \eqn{7.16}, we rewrite it as

\begin{equation}
\qBin{2N+2M+1}{2N}{\sqrt q}=
\sum_{i=0}^{2N}q^{\frac{i}{2}} \qBin{i+M}{i}{q}\qBin{M+2N-i}{M}{q}.
\label{eq:7.19}
\end{equation}
If we let $M\rightarrow\infty$ in \eqn{7.19}, we get a well-known formula for 
the Rogers--Szerg\"o polynomials (Ex.~5, p.~49 in \cite{An5})

\begin{equation}
(-\sqrt q;\sqrt q)_{2N}=
\sum_{i=0}^{2N}q^{\frac{i}{2}} \qBin{2N}{i}{q}.
\label{eq:7.20}
\end{equation}
On the other hand, if we let $N\rightarrow\infty$ in \eqn{7.19}, we find a special case 
of $q$-binomial theorem ((3.3.7) in \cite{An5})

\begin{equation}
\frac{1}{(\sqrt q;q)_{1+M}}=
\sum_{i\geq 0}q^{\frac{i}{2}} \qBin{i+M}{i}{q}.
\label{eq:7.21}
\end{equation}
Thus, formula \eqn{7.19} connects two known, but previously unrelated polynomial identities. 

Next, we apply Burge transform \eqn{7.4} to \eqn{7.16}.
Employing the $q$--Pfaff--Saalsch\"utz sum, we deduce that

\begin{equation}
\qBin{2N+2M}{2N}{\sqrt q}=
\sum_{j=-M}^{M} q^{j^2+\left(\frac{j}{2}\right)} Q(N,M,j,j,q).
\label{eq:7.22}
\end{equation}
If we iterate this last identity by \eqn{7.3}, we find

\begin{equation}
\sum_{i=0}^{M} q^{i^2} \qBin{2N+M-i}{2N}{q} \qBin{2N}{2i}{\sqrt q}=
\sum_{j=-M}^{M} q^{2j^2+\left(\frac{j}{2}\right)} Q(N,M,2j,j,q).
\label{eq:7.23}
\end{equation}
Surprisingly, under the substitution $q\rightarrow q^2$, \eqn{7.23} becomes \eqn{5.6}. 
If we continue to iterate \eqn{7.23} by \eqn{7.3}, we obtain after $\nu-1$ iterations:
$$
\sum_{n_1,\ldots,n_\nu\geq 0} 
q^{N_1^2+\cdots+N_\nu^2} \qBin{2N+M-N_1}{M-N_1}{q}
\prod_{j=1}^\nu \qBin{(1+\delta_{j,\nu})n_j+2N-2\sum_{l=1}^jN_l}
                     {(1+\delta_{j,\nu})n_j}{q^{\frac{2-\delta_{j,\nu}}{2}}}
$$
\begin{equation}
=\sum_{j=-M}^{M} q^{(\nu+1)j^2+\left(\frac{j}{2}\right)}
Q(N,M,(\nu+1)j,j,q),
\label{eq:7.24}
\end{equation}
where, as ususal,
\begin{equation}
N_j=n_j+n_{j+1}+\cdots+n_\nu,
\label{eq:7.25}
\end{equation}
and it is understood that $N\geq N_1+\cdots+N_\nu$.
It is easy to check that \eqn{7.24} remains invariant under transformation 
$q\rightarrow\frac{1}{q}$.

Finally, replacing $q$ by $q^2$ and $\nu$ by $\nu-1$ in \eqn{7.24} and letting $N$ and 
$M$ tend to infinity, we deduce 
with the aid of the Jacobi's triple product formula ((II.28) in \cite{GR}) that
$$
\sum_{n_1,\ldots,n_{\nu-1}\geq 0} 
\frac{q^{2(N_1^2+\cdots+N_{\nu-1}^2)}}
     {(q^2;q^2)_{n_1}\ldots(q^2;q^2)_{n_{\nu-2}}(q)_{2n_{\nu-1}}}=
\frac{1}{(q^2;q^2)_\infty} \sum_{j=-\infty}^\infty q^{2\nu j^2+j}
$$
\begin{equation}
=\frac{1}{(q^{2\nu-1},q^{2\nu+1};q^{4\nu})_\infty}
\frac{(q^{8\nu},q^{4\nu},q^{4\nu-2},q^{4\nu+2};q^{8\nu})_\infty}{(q^2;q^2)_\infty}.
\label{eq:7.26}
\end{equation}
We remark that \eqn{7.26} is the analytic version of the partition theorem in \cite{AS}. 

A complete Burge tree with the ``root" \eqn{7.16} will be described elsewhere.

\medskip

\section{Further observations}
\label{sec:8}

Let us recall the surprise encountered in the last section. Identities \eqn{5.1} and 
\eqn{5.6} can be derived either via the WP--Bailey tree or, alternatively, via the Burge tree. 
To paraphrase, a single iteration of $(A_N^{(1)},B_N^{(1)})$, $(A_N^{(2)},B_N^{(2)})$ 
and $(A_N^{(3)},B_N^{(3)})$  by \eqn{1.10} and \eqn{1.11} gives a result, which is 
equivalent to a double iteration of these pairs by the Burge transform \eqn{7.4}. Clearly, 
one wants a deeper understanding of relation between these two constructs.
To this end we use a $q$--Pfaff-Saalsch\"utz sum to find a following variation on the 
construct \eqn{1.10}, \eqn{1.11}.

\begin{thm}
If $A_N(q)$ and $B_N(M,q)$ form a Burge pair satisfying \eqn{7.1}, then the identity 
$$
\sum_{i\geq 0} q^{(1+2M)i} \qBin{2M+N-i}{N-i}{q}
B_i(-1-M,q)=
$$
\begin{equation}
\sum_{j\ge 0} A_j(q) q^{2aj^2(a-b)} Q(N,M,aj,(2a-b)j,q)
\label{eq:8.1}
\end{equation}
holds true.
\label{thm:2}
\end{thm}

Analogous to what we observed in the case of \eqn{1.10} and \eqn{1.11}, we see that a 
double application of the last theorem gives back an original Burge pair.
However, in case $2a\ge b\ge 0$, one can exploit the symmetry
\begin{equation}
Q(N,M,aj,bj,q)=Q(M,N,bj,aj,q),
\label{eq:8.2}
\end{equation}
to iterate \eqn{7.1} by \eqn{8.1} further in order to get a multiseries representation for 
$$
\sum_{j\ge 0} A_j(q) q^{j^2(a-b)\nu(\nu(a-b)+a+b)} 
Q(N,M,f_{a,b,\nu+1}j,f_{a,b,\nu}j,q),
$$
with $\nu$ being a number of iterations and $f_{a,b,\nu}=\nu(a-b)+b$.
It is worth mentioning that unlike what happens if one iterates by \eqn{7.3} and \eqn{7.4}, 
iteration of \eqn{7.1} by \eqn{8.1} together with \eqn{8.2} results in a one dimensional 
chain of identities, instead of a binary tree.
For example, a chain with the ``seed" \eqn{7.16} gives a self dual identity
$$
\sum_{m_1,\ldots,m_\nu\geq 0} q^{2\sum_{i=1}^{\nu-1}m_i^2+m_\nu^2-2\sum_{i=1}^{\nu-1}m_im_{i+1}}
\qBin{2N+M-m_1}{M-m_1}{q}
$$
\begin{equation}
\times\prod_{i=1}^\nu \qBin{m_{i-1}+m_{i+1}}{2m_i}{q^{\frac{2-\delta_{i,\nu}}{2}}}=
\sum_{j=-\lfloor\frac{M}{\nu}\rfloor} 
^{\lfloor\frac{M}{\nu}\rfloor} 
q^{j^2\nu(\nu+1)+\frac{j}{2}}
Q(N,M,(\nu+1)j,\nu j,q),
\label{eq:8.3}
\end{equation}
with $m_0:=N$ and $m_{\nu+1}:=m_{\nu-1}$.

To emphasize a close relation between our chain and a Burge tree, we observe that the above 
identity could have been derived by applying Burge's transform \eqn{7.4} to \eqn{7.24} with 
$\nu\rightarrow\nu-1$ and by changing the summation variables in the resulting identity.

Despite such strong similarities, \eqn{8.1} and \eqn{7.5}, \eqn{7.6} are clearly different. 
In particular, while a double application of \eqn{8.1} to \eqn{7.1} gives back the initial 
identity, this is not the case for either \eqn{7.5} or \eqn{7.6}. Further implications of the 
Theorem~\ref{thm:2} will be explored in our subsequent publications.

\medskip

\section{Concluding Remarks}
\label{sec:9}

Our main object has been to illustrate the tremendous possibilities nascent in WP--Bailey 
pairs and their iteration via \eqn{1.7}--\eqn{1.11} and \eqn{8.1}, \eqn{8.2}. In a sense, 
the original methods developed by Bailey and employed so successfully by Gasper and Rahman 
(\cite{GR}, \S 2.8), could be used to prove many results in this paper. However the power 
of WP--Bailey pairs lies in the fact that it allows one to see past the leaves to the tree.
Given the large number of new $q$-hypergeometric transformations derived in this paper, 
it is clear that the WP--Bailey tree contains many more applications of interest. 

\subsection*{Acknowledgment}

We would like to thank K.~Alladi, G.~Gasper and S.O.~Warnaar for their interest and 
comments on this manuscript.
Research of G.E.~Andrews was partially supported by NSF Grant DMS-9206993 and
research of A.~Berkovich was supported in part by NSF Grant DMS-0088975.

\end{document}